\input amstex
\documentstyle{amsppt}
%
%
\nopagenumbers
\def\HRes{\operatorname{HRes}}
\def\Cl{\operatorname{Cl}}
\def\compos{\,\raise 1pt\hbox{$\sssize\circ$} \,}
\def\const{\operatorname{const}}
\def\negskp{\hskip -2pt}
\pagewidth{360pt}
\pageheight{606pt}
\topmatter
\title On rational extension of Heisenberg algebra.
\endtitle
\author
Ruslan A. Sharipov
\endauthor
\abstract
Construction of rational extension for Heisenberg algebra with one
pair of generators is discussed.
\endabstract
\address Rabochaya~street~5, 450003, Ufa, Russia
\endaddress
\email \vtop to 20pt{\hsize=280pt\noindent
R\_\hskip 1pt Sharipov\@ic.bashedu.ru\newline
ruslan-sharipov\@usa.net\vss}
\endemail
\urladdr
http:/\negskp/www.geocities.com/CapeCanaveral/Lab/5341
\endurladdr
\endtopmatter
\document
\head
1. Introduction.
\endhead
    Heisenberg algebra with one pair of generators arises in describing
one-dimensio\-nal dynamics of spinless particle in non-relativistic quantum
mechanics. Stationary states of such particle are described by scalar
functions $\Psi(x)$ from the space $L_2(\Bbb R)$:
$$
\int\limits^{+\infty}_{-\infty}|\Psi(x)|^2\,dx<\infty.
$$
Physical parameters of a particle in quantum mechanics are represented
by corresponding operators that can be applied to wave-functions $\Psi$.
Thus coordinate of particle is represented by the operator of multiplication
by $x$, while momentum of particle is represented by differential operator
$$
p=-i\,\hbar\,\frac{d}{dx}.
$$
Operators $x$ and $p$ are bound by commutational relationship $[x,\,p]=-i
\,\hbar$. Denote by $y$ the operator of differentiation with respect to
$x$. This operator differs from the operator of momentum only by constant
factor. For operators $x$ and $y$ we have the following commutational
relationship:
$$
x\,y-y\,x=[x,\,y]=-1.\hskip -2em
\tag1.1
$$
Transferring from quantum mechanics to purely algebraic situation, we
consider associative (but not commutative) algebra with unity $H=H(x,y)$
defined by two generators $x$ and $y$, which are bound by commutational
relationship \thetag{1.1}. It is called {\bf Heisenberg algebra}. In
such treatment $y$ is no longer an operator of differentiation with
respect to $x$, but it is simply a symbol for denoting one of the
generators, playing the same role as a symbol $x$ denoting another
generator.\par
     As a linear space, algebra $H(x,y)$ is generated by its unit element
and by monomials $x^{\alpha_1}\cdot y^{\beta_1}\cdot\ldots\cdot x^{\alpha_n}
\cdot y^{\beta_n}$ with entire non-negative powers. Due to the commutational
relationship \thetag{1.1} such monomial can be transformed to the linear
combination of monomials of the form $x^{\alpha}\cdot y^{\beta}$. Therefore,
as a linear space, algebra $H(x,y)$ can be identified with the set of
polynomials of two variables $\Bbb C[x,y]$. However, the operation of
multiplication in $H(x,y)$ differs from multiplication of polynomials
in $\Bbb C[x,y]$. Let's study this operation in more details. First
let's calculate commutator $[x,\,y^q]$, where power $q$ is some entire
non-negative number. For $q\geqslant 2$ we get
$$
[x,\,y^q]=[x,\,y\cdot\ldots\cdot y]=[x,\,y]\cdot y^{q-1}+
y\cdot[x,\,y]\cdot y^{q-2}+\ldots+y^{q-1}\cdot[x,\,y]=-q\,y^{q-1}.
$$
Denote $f(y)=y^q$, then $q\,y^{q-1}=f'(y)$. Now the above formula
is rewritten as
$$
[x,\,f(y)]=-f'(y).\hskip -2em
\tag1.2
$$
Taking for $f$ the polynomial of two variables $f(x,y)$, one can generalize
formula \thetag{1.2}. Namely we have the following two formulas:
$$
\xalignat 2
&[x,\,f]=-f'_y,&&[y,\,f]=f'_x.\hskip -2em
\tag1.3
\endxalignat
$$
Here, when substituting generators of Heisenberg algebra $x$ and $y$
into polynomial $f(x,y)$, they are assumed to be ordered in a following
way:
$$
f(x,y)=\sum_{p\,q}f_{pq}\ x^p\cdot y^q.\hskip -2em
\tag1.4
$$
Derivation of formulas \thetag{1.3} is quite similar to derivation of
previous formulas \thetag{1.2}.\par
    In formulas \thetag{1.3} we see that transposition of the element
$f=f(x,y)$ with generators $x$ and $y$ in Heisenberg algebra leads to
the differentiation of corresponding polynomial $f(x,y)$. Formulas 
\thetag{1.3} admit further generalization:
$$
\aligned
&f\cdot x^k=\sum^k_{\alpha=0}C^\alpha_k\ 
x^\alpha\cdot D_y^{k-\alpha} f,\\
&y^k\cdot f=\sum^k_{\alpha=0}C^\alpha_k\ D_x^{k-\alpha} f
\cdot y^\alpha.
\endaligned\hskip -2em
\tag1.5
$$
Here $D_x$ and $D_y$ are operators of differentiation in $x$ and $y$
respectively, while $C^\alpha_k$ in formulas \thetag{1.5} are binomial
coefficients:
$$
C^\alpha_k=\frac{k!}{\alpha!\ (k-\alpha)!}.\hskip -2em
\tag1.6
$$
Formulas \thetag{1.5} are proved by induction in $k$. For $k=0$ they
are obvious. For $k=1$ they reduces to \thetag{1.3}. Inductive step
from $k$ to $k+1$ is provided by the following well-known identity for
binomial coefficients:
$$
C^\alpha_{k+1}=\cases C^{\alpha-1}_k+C^\alpha_k &\text{\ \ for \ }
0<\alpha<k+1,\\
C^{\alpha-1}_k&\text{\ \ for \ }\alpha=k+1,\\
C^\alpha_k&\text{\ \ for \ }\alpha=0.
\endcases
$$
Let's substitute $f=y^{k+q}$ and $f=x^{k+q}$, where $q\geqslant 0$,
into \thetag{1.5}. This yields
$$
\align
&y^{k+q}\cdot x^k=\sum^k_{\alpha=0}
\frac{k!\ (k+q)!}{\alpha!\ (k-\alpha)!\
(q+\alpha)!}\ x^\alpha\cdot y^{q+\alpha},
\hskip -2em
\tag1.7\\
&y^k\cdot x^{k+q}=\sum^k_{\alpha=0}
\frac{k!\ (k+q)!}{\alpha!\ (k-\alpha)!\
(q+\alpha)!}\ x^{\alpha+q}\cdot y^\alpha.
\hskip -2em
\tag1.8
\endalign
$$
Formulas \thetag{1.7} and \thetag{1.8} can be united into one formula
if we write them as follows:
$$
y^q\cdot x^r=\sum\Sb\sssize 0\leqslant \alpha\leqslant q\\
\vspace{1pt}
\sssize 0\leqslant \alpha\leqslant r\endSb
C^\alpha_{qr}\ x^{q-\alpha}\cdot y^{r-\alpha}.
\hskip -2em
\tag1.9
$$
Coefficients $C^\alpha_{qr}$ in this relationship \thetag{1.9} are
determined by formula
$$
C^\alpha_{qr}=\frac{q!\ r!}{(q-\alpha)!\ \alpha!\ (r-\alpha)!}.
\hskip -2em
\tag1.10
$$
Formula \thetag{1.10} for $C^\alpha_{qr}$ is similar to formula
\thetag{1.6} for binomial coefficients. This is why we have
chosen symbols $C^\alpha_{qr}$ for coefficients in \thetag{1.9}.
\head
2. Some properties of Heisenberg algebra.
\endhead
    Consider two polynomials of the form \thetag{1.4}. Suppose that
these are $f$ and $g$:
$$
\xalignat 2
&\quad f=\sum^n_{p=0}\sum^n_{q=0}f_{pq}\ x^p\cdot y^q,
&&g=\sum^n_{r=0}\sum^n_{s=0}g_{rs}\ x^r\cdot y^s.
\hskip -2em
\tag2.1
\endxalignat
$$
Let's calculate the product of polynomials \thetag{2.1} in Heisenberg
algebra:
$$
f\cdot g=\sum^n_{p=0}\sum^n_{q=0}f_{pq}\ x^p\cdot y^q\cdot g(x,y)=
\sum^n_{p=0}\sum^n_{q=0}\sum^q_{\alpha=0}f_{pq}\ C^\alpha_q\ x^p
\cdot D^\alpha_x g(x,y)\cdot y^{q-\alpha}.
$$
Here we used second relationship \thetag{1.5} and took into account
the following symmetry of binomial coefficients: $C^\alpha_q
=C^{q-\alpha}_q$. Further let's note that 
$$
\frac{D^\alpha_y(y^q)}{\alpha!}
=\cases C^\alpha_q\ y^{q-\alpha}
&\text{\ \ for \ }\alpha\leqslant q,\\
\vspace{2ex}
\ \ \ 0&\text{\ \ for \ }\alpha>q.
\endcases
$$
Therefore if we construct by $f$ and $g$ a new polynomial
$$
h(x,y)=\sum^{n+1}_{\alpha=0}\frac{D^\alpha_y f\ D^\alpha_x g}{\alpha!},
$$
and if, upon collecting similar terms in it, we arrange variables $x$ and
$y$ in a natural order as in \thetag{1.4}, then for the product of
polynomials $f\cdot g$ we can write $f\cdot g=h$:
$$
\pagebreak
f\cdot g=\sum^{n+1}_{\alpha=0}\frac{D^\alpha_y f\ D^\alpha_x g}{\alpha!}.
\hskip -2em
\tag2.2
$$
For us it's important that formula \thetag{2.2} reduces multiplication of
polynomials in Heisenberg algebra to the ordinary multiplication of
polynomials in the ring $\Bbb C[x,y]$ and to the differentiation of these
polynomials. If we do not restrict the degree of polynomials $f$ and $g$
by particular number $n$, we can rewrite formula \thetag{2.2} as
$$
f\cdot g=\sum^\infty_{\alpha=0}\frac{D^\alpha_y f\ D^\alpha_x g}{\alpha!}.
\hskip -2em
\tag2.3
$$
\proclaim{Theorem 2.1} Heisenberg algebra with one pair of generators
$H(x,y)$ has no divisors of zero.
\endproclaim
\demo{Proof}  Recall that two nonzero elements $f$ and $g$ are called
divisors of zero if their product is zero: $f\cdot g=0$ (see \cite{1}
or \cite{2}). Suppose that such two elements in $H(x,y)$ do exist.
They should be represented by two polynomials
$$
\xalignat 2
&\quad f=\sum^m_{q=0}f_q(x)\cdot y^q,
&&g=\sum^n_{q=0}g_q(x)\cdot y^q,
\hskip -2em
\tag2.4
\endxalignat
$$
where $f_m(x)\neq 0$ and $g_n(x)\neq 0$. Substituting polynomials
\thetag{2.4} into the formula \thetag{2.3}, for the product of
these polynomials we get
$$
h=\sum^{m+n}_{q=0}h_q(x)\cdot y^q.\hskip -2em
\tag2.5
$$
For the leading term in polynomial \thetag{2.5} we have $h_{m+n}(x)=f_m(x)
\ g_n(x)$. Therefore from $f_m(x)\neq 0$ and $g_n(x)\neq 0$ we derive
$h_{m+n}(x)\neq 0$. Hence $h\neq 0$. This contradicts to the assumption
that $f\cdot g=0$. Theorem is proved.\qed\enddemo
    If the equalities $g_1\cdot f=1$ and $f\cdot g_2=1$ are fulfilled,
then element $g_1$ is called {\bf left inverse element} for $f$, and
$g_2$ is called {\bf right inverse element} for $f$. It's easy to show
that if element $f$ in associative algebra has both left and right
inverse elements, then these two inverse elements do coincide (see
\cite{1} or \cite{2}):
$$
g_1=g_1\cdot (f\cdot g_2)=(g_1\cdot f)\cdot g_2=g_2.
$$
In associative algebra without divisors of zero the existence of
left inverse element $g_1$ for $f$ implies the existence of right
inverse element $g_2=g_1$. Indeed, from $1=g_1\cdot f$ it follows
that $f=f\cdot (g_1\cdot f)=(f\cdot g_1)\cdot f$. Then we have
$(f\cdot g_1-1)\cdot f=0$. The element $f\neq 0$ cannot be a divisor
of zero, therefore $f\cdot g_1-1=0$. This yields the relationship
$f\cdot g_1=1$.\par
    And conversely, the existence of right inverse element $g_2$ for $f$
implies the existence of left inverse element $g_1=g_2$. This fact is
proved similarly. From $1=f\cdot g_2$ it follows that $f=(f\cdot g_2)\cdot
f=f\cdot (g_2\cdot f)$. Further we have $f\cdot(1-g_2\cdot f)=0$. This
yields the required relationship $g_2\cdot f=1$.\par
    If element $g$ is both left and right inverse for the element $f$, then
$g$ is called {\bf bilateral inverse element} for $f$ or simply {\bf inverse
element}. It is denoted as $g=f^{-1}$. Element $f$ possessing bilateral
inverse element $f^{-1}$ is called {\bf invertible}.\par
\proclaim{Theorem 2.2} Element $f$ of Heisenberg algebra $H(x,y)$ is
invertible if and only if it is in the field of scalars.
\endproclaim
\demo{Proof} Suppose that $f$ is invertible and $f\notin\Bbb C$. Then it
is represented by a polynomial $f(x,y)$ which is not constant. Suppose,
for instance, that $f(x,y)$ has an actual entry of the variable $y$. Then
$$
f=\sum^m_{q=0}f_q(x)\cdot y^q,\hskip -2em
\tag2.6
$$
where $m\neq 0$ and $f_m(x)\neq 0$. Let $g=f^{-1}$. We write polynomial
$g$ as 
$$
g=\sum^n_{q=0}g_q(x)\cdot y^q,\hskip -2em
\tag2.7
$$
where $g_n(x)\neq 0$. Substituting polynomials \thetag{2.6} and \thetag{2.7}
into formula \thetag{2.3}, for the product of these polynomials $h=f\cdot g$
we get formula \thetag{2.5} with $n+m\neq 0$ and $h_{n+m}(x)=f_m(x)\ g_n(x)$.
Hence $h_{n+m}(x)\neq 0$. This result contradicts to the equality
$h=f\cdot g=f\cdot f^{-1}=1$.\par
    Now we have to consider the case, when $f(x,y)$ doesn't contain actual
entries of the variable $y$. Then $f=f(x)\notin\Bbb C$. In this case the
product of polynomials $f\cdot g$ in Heisenberg algebra coincides with their
product in the ring $\Bbb C[x,y]$:
$$
f\cdot g=f\ g=1.
$$
But non-constant polynomial $f(x)$ isn't invertible element of the ring
$\Bbb C[x,y]$. The contradictions obtained prove that $f=\const\in\Bbb C$.
\qed\enddemo
\definition{Definition 2.1} Two nonzero elements $f$ and $g$ in
noncommutative algebra are called {\bf left comeasurable} if one can
find nonzero elements $a$ and $b$ such that $a\cdot f=b\cdot g$.
\enddefinition
\definition{Definition 2.2} Two nonzero elements $f$ and $g$ in
noncommutative algebra are called {\bf right comeasurable} if one can
find nonzero elements $a$ and $b$ such that $f\cdot a=g\cdot b$.
\enddefinition
    In commutative algebra the concept of comeasurability is trivial,
since any two elements in commutative algebra are comeasurable. In this
aspect noncommutative Heisenberg algebra is similar to commutative
ones. Namely we have the theorem.
\proclaim{Theorem 2.3} Arbitrary two nonzero elements $f$ and $g$ in
Heisenberg algebra $H(x,y)$ are left comeasurable.
\endproclaim
\demo{Proof} Elements of Heisenberg algebra are represented by polynomials.
Their product is determined by formula \thetag{2.3}. Suppose that elements
$f$ and $g$ are given by polynomials \thetag{2.6} and \thetag{2.7}, where
$m$ and $n$ are degrees of these polynomials with respect to the variable
$y$. For elements $a$ and $b$ in definition~2.1 (if they exist) from the
equality $a\cdot f=b\cdot g$ we derive
$$
m+\deg_y(a)=n+\deg_y(b).
$$
Therefore we shall construct the elements $a$ and $b$ in form of
polynomials
$$
\xalignat 2
&\quad a=\sum^n_{k=0}a_k(x)\cdot y^k,
&&b=\sum^m_{k=0}b_k(x)\cdot y^k.
\hskip -2em
\tag2.8
\endxalignat
$$
Applying formula \thetag{2.3}, we transform the equality $a\cdot f
=b\cdot g$, which should be satisfied by polynomials \thetag{2.6},
\thetag{2.7}, and \thetag{2.8}:
$$
\sum^n_{k=0}\frac{D^k_y a\ D^k_x f}{k!}=\sum^m_{k=0}\frac{D^k_y b\
D^k_x g}{k!}.\hskip -2em
\tag2.9
$$
Now \thetag{2.9} is an ordinary polynomial equality in the ring
$\Bbb C[x,y]$. Let's do further transformations in this equality:
$$
\sum^n_{k=0}\sum^n_{q=k}\sum^m_{s=0}C^k_q\ a_q(x)\ D^k_xf_s(x)\ y^{q-k+s}
=\sum^m_{k=0}\sum^m_{q=k}\sum^n_{s=0}C^k_q\ b_q(x)\ D^k_xg_s(x)\ y^{q-k+s}.
$$
First we change order of summation and replace index $k$ by an index
$p=q-k$:
$$
\sum^m_{s=0}\sum^n_{q=0}\sum^q_{p=0}C^p_q\ a_q(x)\ D^{q-p}_xf_s(x)\ y^{p+s}
=\sum^n_{s=0}\sum^m_{q=0}\sum^q_{p=0}C^p_q\ b_q(x)\ D^{q-p}_xg_s(x)\ y^{p+s}.
$$
Then we do another change of the order of summation and replace $s$ by
$r=p+s$:
$$
\aligned
\sum^{m+n}_{r=0}&\left(\ \sum^{\min(n,r)}_{p=\max(0,r-m)}\left(\,
\sum^n_{q=p}C^p_q\ a_q(x)\ D^{q-p}_xf_{r-p}(x)\right)\right)y^r=\\
&=\sum^{m+n}_{r=0}\left(\ \sum^{\min(m,r)}_{p=\max(0,r-n)}\left(\,
\sum^m_{q=p}C^p_q\ b_q(x)\ D^{q-p}_xg_{r-p}(x)\right)\right)y^r.
\endaligned\hskip -2em
\tag2.10
$$
Here \thetag{2.10} is an equality of two polynomials of the order
$m+n$ with respect to variable $y$. In their leading terms with
$r=n+m$ we find
$$
a_n(x)\,f_m(x)=b_m(x)\,g_n(x).\hskip -2em
\tag2.11
$$
The equation \thetag{2.11} can be easily satisfied if we define
$a_n$ and $b_m$ by formulas
$$
\xalignat 2
&\quad a_n(x)=g_n(x)\,\varphi(x),&&b_m(x)=f_m(x)\,\varphi(x).\hskip -2em
\tag2.12
\endxalignat
$$
We shall determine polynomial $\varphi(x)$ later. Now let's consider
again the polynomial equality \thetag{2.10}. In leading order $r=n+m$
it is fulfilled due to \thetag{2.12}. Equating coefficients of $y^r$
for other powers, we get $n+m$ equalities which can be treated as the
equations with respect to $n+m$ polynomials $a_0(x),\,\ldots,\,
a_{m-1}(x)$ and $b_0(x),\,\ldots,\,b_{n-1}(x)$. Note that \thetag{2.10}
doesn't contain the derivatives of polynomials $a_0(x),\,\ldots,\,a_{m-1}
(x)$ and $b_0(x),\,\ldots,\,b_{n-1}(x)$. Therefore we have the
non-homogeneous system of linear algebraic equations with respect to
these polynomials. It can be written in matrix form as follows
$$
\Vmatrix M^1_1&\hdots&M^1_n&M^1_{n+1}&\hdots&M^1_{n+m}\\
\vspace{1.5ex}
\vdots&\ddots&\vdots&\vdots&\ddots&\vdots\\
\vspace{1.5ex}
M^n_1&\hdots&M^n_n&M^n_{n+1}&\hdots&M^n_{n+m}\\
\vspace{1.5ex}
M^{n+1}_1&\hdots&M^{n+1}_n&M^{n+1}_{n+1}&\hdots&M^{n+1}_{n+m}\\
\vspace{1.5ex}
\vdots&\ddots&\vdots&\vdots&\ddots&\vdots\\
\vspace{1.5ex}
M^{n+m}_1&\hdots&M^{n+m}_n&M^{n+m}_{n+1}&\hdots&M^{n+m}_{n+m}
\endVmatrix\
\Vmatrix a_{n-1}\\ \vspace{1.5ex}\vdots\\ \vspace{1.5ex}a_0\\
\vspace{1.5ex}b_{m-1}\\ \vspace{1.5ex}\vdots\\ \vspace{1.5ex}b_0
\endVmatrix=
\Vmatrix A_{n-1}\\ \vspace{1.5ex}\vdots\\ \vspace{1.5ex}A_0\\
\vspace{1.5ex}B_{m-1}\\ \vspace{1.5ex}\vdots\\ \vspace{1.5ex} B_0
\endVmatrix.\hskip -2em
\tag2.13
$$
Elements of matrix $M$ in \thetag{2.13} are polynomials with respect
to the variable $x$, they are determined by coefficients $f_q(x)$ and
$g_q(x)$ in polynomials \thetag{2.6} and \thetag{2.7}. Quantities
$A_0,\,\ldots,\,A_{m-1}$ and $B_0,\,\ldots,\,B_{n-1}$ in right hand
side of \thetag{2.13} are also polynomials in $x$. They are determined
by coefficients $f_q(x)$ and $g_q(x)$ in polynomials \thetag{2.6} and
\thetag{2.7}, and by our choice of leading coefficients $a_n(x)$ and
$b_m(x)$ in \thetag{2.8}. In particular, if we determine $a_n(x)$ and
$b_m(x)$ by formula \thetag{2.12}, then polynomials $A_0,\,\ldots,\,
A_{m-1}$ and $B_0,\,\ldots,\,B_{n-1}$ has common factor $\varphi(x)$.
\par
    First let's study the case, when $\mu(x)=\det M\neq 0$. The quantity
$\mu(x)$ can be treated as {\bf noncommutative analog of resultant} in
Heisenberg algebra for two polynomials $f(x,y)$ and $g(x,y)$ respective
to the variable $y$:
$$
\mu(x)=\HRes_y(f,g).
$$
For the case $\mu(x)\neq 0$ we choose $\varphi(x)=\mu(x)$ in \thetag{2.12}. 
In this case system of equations \thetag{2.13} has unique solution given
by polynomials in $x$. So the required polynomials \thetag{2.8} are
constructed.\par
    Now let's consider degenerate case $\mu(x)=\det M=0$. In this case
we choose $\varphi(x)=0$ in \thetag{2.12}. This choice makes zero the
right hand sides of the equations \thetag{2.13}. These equations become
homogeneous. Homogeneous system of linear algebraic equations with
degenerate square matrix always has at least one nontrivial solution.
We can find it applying Gauss's method. In general case this solution
$a_0(x),\,\ldots,\,a_{m-1}(x),\,b_0(x),\,\ldots,\,b_{n-1}(x)$ will be
represented by rational functions. However, the solution of homogeneous
system of linear equations is always determined only up to an arbitrary
common factor. We can choose this factor so that $a_0(x),\,\ldots,\,
a_{m-1}(x),\,b_0(x),\,\ldots,\,b_{n-1}(x)$ will become polynomials.
Thus, the required polynomials \thetag{2.8} do exist in degenerate case
$\HRes_y(f,g)=\det M=0$ as well. Theorem is proved.
\qed\enddemo
\head
3. Rational extension of Heisenberg algebra.
\endhead
    It is known that rational functions are determined as ratio of
two polynomials. In abstract algebraic situation for commutative rings
and algebras this is generalized in form of the construction of the
field of fractions (see \cite{1} and \cite{2}). In the case of Heisenberg
algebra we need noncommutative analog of this construction.\par
    Consider a set $\Cal M$ consisting of ordered pairs of elements from
Heisenberg algebra. We shall denote such pairs as $b^{-1}\compos a$, where
$a$ and $b$ are two elements of $H(x,y)$. Into $\Cal M$ we include only
those pairs $b^{-1}\compos a$ for  which $b\neq 0$. Two pairs $b^{-1}
\compos a$ and $d^{-1}\compos c$ are called {\bf equivalent}
if there are two elements $u$ and $v$ in $H(x,y)$ such that
$$
\xalignat 2
&u\cdot b=v\cdot d,&&u\cdot a=v\cdot c.\hskip -2em
\tag3.1
\endxalignat
$$
Properties of reflexivity and symmetry for the equivalence relation just
introduced are obvious. Let's check the property of transitivity. Suppose
that $b^{-1}\compos a\sim d^{-1}\compos c$ and $d^{-1}\compos c\sim f^{-1}
\compos e$. Then there are some elements $u$, $v$, $w$ and $z$ in $H(x,y)$
such that the following relationships are fulfilled:
$$
\xalignat 2
&u\cdot b=v\cdot d,&&u\cdot a=v\cdot c,\\
&w\cdot d=z\cdot f,&&w\cdot c=z\cdot e,
\endxalignat
$$
Elements $v$ and $w$ are nonzero. This follows from $b\neq 0$, $d\neq 0$,
and $f\neq 0$ due to the theorem~2.1. Let's apply theorem~2.3 to $v$ and
$w$. It asserts that there are two elements $p$ and $q$ in $H(x,y)$ such
that $p\cdot v=q\cdot w$. Therefore we can write
$$
\gather
(p\cdot u)\cdot b=p\cdot v\cdot d=q\cdot w\cdot d=(q\cdot z)\cdot f,\\
(p\cdot u)\cdot a=p\cdot v\cdot c=q\cdot w\cdot c=(q\cdot z)\cdot e.
\endgather
$$
From the above relationships we get $b^{-1}\compos a\sim f^{-1}\compos e$.
This means that the property transitivity is present.\par
    Thus the set $\Cal M$ consisting of pairs $b^{-1}\compos a$, which
further will be called fractions, breaks into classes of equivalence.
Denote by $\Cal H(x,y)$ the set of such classes. The class of equivalence
for the fraction $b^{-1}\compos a$ might be denoted by $\Cl(b^{-1}\compos
a)$. However, for sake of simplicity the sign of class is usually
omitted. The sign of equivalence $\sim$ also is often replaced by the
sign of equality, emphasizing that equivalent fractions are treated as
unseparable objects.\par
\proclaim{Lemma 3.1} If $b^{-1}\compos a\sim d^{-1}\compos c$, then from
the equality of denominators follows the equality of numerators $a=c$.
\endproclaim
    Indeed, equivalence $b^{-1}\compos a\sim d^{-1}\compos c$ means that
there exist two elements $u$ and $v$ such that the relationships
\thetag{3.1} are fulfilled. Due to $b=d\neq 0$ first of them can be
brought to the form $(u-v)\cdot b=0$. According to the theorem~2.1,
Heisenberg algebra has no divisors of zero. Therefore $u=v\neq 0$. Now
second relationship \thetag{3.1} yields $u(a-c)=0$. Hence $a=c$.\par
    Let's define algebraic operations in factor set $\Cal H(x,y)=\Cal M
/\kern -4pt\sim$. Let's begin with the operation of summation. Suppose
that we have two fractions. Their denominators are nonzero. We apply
theorem~2.3 to them and find elements $u$ and $v$ from $H(x,y)$ such
that the following equality is fulfilled:
$$
u\cdot b=v\cdot d=f\neq 0.\hskip -2em
\tag3.2
$$
Now we can bring fractions $b^{-1}\compos a$ and $d^{-1}\compos c$ to
the common denominator:
$$
\pagebreak
\xalignat 2
&b^{-1}\compos a=(u\cdot b)^{-1}\compos (u\cdot a),
&&d^{-1}\compos c=(v\cdot d)^{-1}\compos (v\cdot c).
\endxalignat
$$
Now we define the sum of two fractions $b^{-1}\compos a$ and
$d^{-1}\compos c$ as
$$
b^{-1}\compos a+d^{-1}\compos c=f^{-1}\compos(u\cdot a+v\cdot c),\hskip -2em
\tag3.3
$$
where $u$, $v$, and denominator $f$ are determined by the relationship
\thetag{3.2}. One should check correctness of determining sum of fractions
by means of formula \thetag{3.3}.
\proclaim{Lemma 3.2} If fraction $b^{-1}\compos a$ is equivalent to
$\tilde b^{-1}\compos\tilde a$, and fraction $d^{-1}\compos c$ is
equivalent to $\tilde d^{-1}\compos\tilde c$, then sum of fractions 
$b^{-1}\compos a+d^{-1}\compos c$ constructed by formula
\thetag{3.3} is equivalent to the sum $\tilde b^{-1}\compos\tilde a
+\tilde d^{-1}\compos\tilde c$ constructed by the same formula.
\endproclaim
    Idea of proof consist in bringing fractions $f^{-1}\compos(u\cdot
a+v\cdot c)$ and $\tilde f^{-1}\compos(\tilde u\cdot\tilde a+\tilde v
\cdot\tilde c)$, arising due to formula \thetag{3.3}, to the common
denominator by applying theorem~2.3. Thereby fractions $b^{-1}\compos a$,
\ $d^{-1}\compos c$, \ $\tilde b^{-1}\compos\tilde a$, \ $\tilde d^{-1}
\compos\tilde c$ are also brought to the same denominator. Further
steps of proof consist in applying lemma~3.1 and in elementary
calculations.\par
    Note that in the set of fractions $\Cal H(x,y)=\Cal M/\kern -4pt\sim$
there is a zero element. It is represented by fractions with numerator
$a=0$ and arbitrary denominator $b\neq 0$.\par
    Now let's determine {\bf multiplication} in the set of fractions
$\Cal H(x,y)$. Suppose that we have two fractions $b^{-1}\compos a$ and
$d^{-1}\compos c$. If $a=0$, then for the product of these fractions
we set by definition
$$
(b^{-1}\compos a)\cdot (d^{-1}\compos c)=0.\hskip -2em
\tag3.4
$$
Suppose $a\neq 0$. Then $d\neq 0$ and $d\cdot a\neq 0$. Let's apply
theorem~2.3 to this pair of nonzero elements. As a result we find
elements $u$ and $v$ such that 
$$
u\cdot d=v\cdot (d\cdot a).\hskip -2em
\tag3.5
$$
From the equality \thetag{3.5} it follows that 
$$
b^{-1}\compos a\sim(v\cdot d\cdot b)^{-1}\compos
(v\cdot d\cdot a)\sim(v\cdot d\cdot b)^{-1}\compos(u\cdot d).
$$
These relationships could be a motivation for determining the
product of fractions $b^{-1}\compos a$ and $d^{-1}\compos c$
by means of formula
$$
(b^{-1}\compos a)\cdot (d^{-1}\compos c)=(v\cdot d\cdot b)^{-1}
\compos(u\cdot c),\hskip -2em
\tag3.6
$$
where elements $u$ and $v$ are determined according to the formula
\thetag{3.5}. One should check correctness of defining the operation
of multiplication by formula \thetag{3.6}. For the formula
\thetag{3.4}, which defines the product of fractions in the case $a=0$,
the checking procedure is trivial.
\proclaim{Lemma 3.3} If fraction $b^{-1}\compos a$ is equivalent
to the fraction $\tilde b^{-1}\compos\tilde a$, and fraction
$d^{-1}\compos c$ is equivalent to the fraction $\tilde d^{-1}\compos
\tilde c$, then the product $(b^{-1}\compos a)\cdot (d^{-1}\compos c)$
constructed by formula \thetag{3.6} is equivalent to the product
$(\tilde b^{-1}\compos\tilde a)\cdot(\tilde d^{-1}\compos\tilde c)$
constructed by the same formula.
\endproclaim
\demo{Proof} The use of formula \thetag{3.6} for multiplying fractions
listed in the statement of lemma~3.3 assumes that we find nonzero elements
$u$, $v$, $\tilde u$, $\tilde v$ such that
$$
\xalignat 2
&u\cdot d=v\cdot d\cdot a,&&\tilde u\cdot\tilde d=\tilde v\cdot\tilde d
\cdot\tilde a.\hskip -2em
\tag3.7
\endxalignat
$$
Their existence is granted by theorem~2.3. Under these conditions we
have
$$
\aligned
&z=(b^{-1}\compos a)\cdot (d^{-1}\compos c)=(v\cdot d\cdot b)^{-1}
\compos(u\cdot c),\\
&\tilde z=(\tilde b^{-1}\compos\tilde a)\cdot (\tilde d^{-1}\compos
\tilde c)=(\tilde v\cdot\tilde d\cdot\tilde b)^{-1}\compos(\tilde u
\cdot\tilde c).
\endaligned\hskip -2em
\tag3.8
$$
Let's prove the equivalence of fractions $z$ and $\tilde z$ in
\thetag{3.8}. In order to prove this fact we shall bring $z$ and
$\tilde z$ to common denominator. However, we shall do it in two
steps. First we apply theorem~2.3 and find elements $p$ and $\tilde p$
such that
$$
p\cdot v\cdot d=\tilde p\cdot\tilde v\cdot\tilde d=w.\hskip -2em
\tag3.9
$$
Multiplying numerators and denominators of fractions $z$ and $\tilde z$
by $p$ and $\tilde p$ respectively, they can be transformed as follows:
$$
\xalignat 2
&z\sim (w\cdot b)^{-1}\compos(p\cdot u\cdot c),
&&\tilde z\sim (w\cdot\tilde b)^{-1}\compos(\tilde p\cdot\tilde u
\cdot\tilde c).\hskip -2em
\tag3.10
\endxalignat
$$
Moreover, we multiply the equalities \thetag{3.7} by $p$ and $\tilde p$
respectively and take into account the above relationship \thetag{3.9}:
$$
\xalignat 2
&p\cdot u\cdot d=w\cdot a,&&\tilde p\cdot\tilde u\cdot\tilde d
=w\cdot\tilde a.\hskip -2em
\tag3.11
\endxalignat
$$
In the second step, applying theorem~2.3 once more, we find $q$ and
$\tilde q$ such that
$$
q\cdot w\cdot b=\tilde q\cdot w\cdot\tilde b=m.\hskip -2em
\tag3.12
$$
This allows us to bring fractions \thetag{3.10} to common denominator:
$$
\xalignat 2
&z\sim m^{-1}\compos(q\cdot p\cdot u\cdot c),
&&\tilde z\sim m^{-1}\compos(\tilde q\cdot \tilde p\cdot\tilde u
\cdot\tilde c).\hskip -2em
\tag3.13
\endxalignat
$$
Now let's note that the equality \thetag{3.12} can be used for
bringing fractions $b^{-1}\compos a$ and $\tilde b^{-1}\compos\tilde a$
to common denominator $m$:
$$
\xalignat 2
&b^{-1}\compos a\sim m^{-1}\compos (q\cdot w\cdot a),
&&\tilde b^{-1}\compos\tilde a\sim m^{-1}\compos (\tilde q\cdot
w\cdot\tilde a).
\endxalignat
$$
We know that fractions $b^{-1}\compos a$ and $\tilde b^{-1}\compos\tilde a$
are equivalent. Using lemma~3.1, we get
$$
q\cdot w\cdot a=\tilde q\cdot w\cdot\tilde a.\hskip -2em
\tag3.14
$$
Let's multiply the equalities \thetag{3.11} from left hand side by $q$ and
$\tilde q$ respectively. Then take into account the equality \thetag{3.14}.
This yields the relationship
$$
q\cdot p\cdot u\cdot d=\tilde q\cdot\tilde p\cdot\tilde u
\cdot\tilde d=h.\hskip -2em
\tag3.15
$$
The relationship \thetag{3.15} can be used for bringing fractions
$d^{-1}\compos c$ and $\tilde d^{-1}\compos\tilde c$ to common
denominator $h$:
$$
\xalignat 2
&d^{-1}\compos c\sim h^{-1}\compos (q\cdot p\cdot u\cdot c),
&&\tilde d^{-1}\compos\tilde c\sim h^{-1}\compos (\tilde q\cdot
\tilde p\cdot\tilde u\cdot\tilde c).
\endxalignat
$$
Let's recall that $d^{-1}\compos c\sim\tilde d^{-1}\compos\tilde c$.
Applying lemma~3.1, we now get
$$
q\cdot p\cdot u\cdot c=\tilde q\cdot\tilde p\cdot\tilde u\cdot\tilde c.
\hskip -2em
\tag3.16
$$
Comparing \thetag{3.16} and \thetag{3.13} we conclude that fractions
$z$ and $\tilde z$ are equivalent. This is the very result we were
to obtain.\qed\enddemo
{\bf Conclusions}. Having defined the operations of addition and
multiplication in factor set $\Cal M/\kern -4pt\sim\kern 4pt=\Cal
H(x,y)$, we turn it into noncommutative ring. Moreover, we can keep
the operation of multiplication by scalars $\alpha\in\Bbb C$:
$$
\alpha\ (b^{-1}\compos a)=b^{-1}\compos (\alpha\,a).
$$
Therefore $\Cal H(x,y)$ is an algebra. It is natural to call it
the {\bf rational extension} of Heisenberg algebra $H(x,y)$. Algebra
$H(x,y)$ is embedded into $\Cal H(x,y)$ in form of fractions with
unitary denominator: $a\mapsto (1)^{-1}\compos a$.\par
   Note that any nonzero element in $\Cal H(x,y)$ is invertible:
nonzero fraction $b^{-1}\compos a$ has an inverse fraction $a^{-1}
\compos b$. Therefore $\Cal H(x,y)$ is an {\bf algebra with division}.
It is known (see. \cite{3}) that finite-dimensional associative
algebras with division over the field of reals $\Bbb R$ are exhausted
by two examples: these are complex numbers $\Bbb C$ and algebra of
quaternions $\Bbb K$. Algebra $\Cal H(x,y)$ is an example of an
algebra with division which is, though infinite-dimensional, but
rather effectively constructed.\par
    Recently I.~Z.~Golubchik said me that the above construction of
algebra $\Cal H(x,y)$ is a special case of more general construction
applicable to so called "skew-polynomial algebras" and to abstract
associative algebras satisfying some rather natural conditions
(see \cite{4} and \cite{5} for more details). Nevertheless, I think
that constructive algorithm for finding the ``comeasuring factors''
and the concept of ``noncommutative resultant'' suggested in the
proof of theorem~2.3 could be worth for actual calculations in
extended Heisenberg algebra $\Cal H(x,y)$ and in developing computer
programs for such calculations.
\head
4. Acknowledgments.
\endhead
    I am grateful to B.~I.~Suleymanov, who let me know some unsolved
problems indicated in papers \cite{6--9}. This paper is aimed to
prepare background for solving some of them. I am also grateful
to I.~Z.~Golubchik, who communicated me references to books \cite{4}
and \cite{5}.\par
     This work is supported by grant from Russian Fund for Basic
Research (project No\nolinebreak\.~00\nolinebreak-01-00068,
coordinator Ya\.~T.~Sultanaev), and by grant from Academy of
Sciences of the Republic Bashkortostan (coordinator N.~M.~Asadullin).
I am grateful to these organizations for financial support.\par

\enddocument
\end